%% file: Main.tex
\documentclass[10pt, final, twocolumn, oneside, conference]{IEEEtran}

\usepackage{times}
\usepackage{amsmath}
\usepackage{amsthm}
\usepackage{mathtools}
\usepackage{epsfig}
\usepackage{graphicx}
\usepackage{enumerate} 
\usepackage{epstopdf}
\usepackage{xspace}
\usepackage{url}
\usepackage{makeidx}
\usepackage[usenames]{color}
\usepackage[dvipsnames]{xcolor}
\usepackage{multirow,makecell}
\usepackage{spreadtab}
\usepackage{color}
\usepackage{subfigure}
\usepackage{adjustbox}
\usepackage{listings}
\usepackage{color}
\usepackage{placeins}
\usepackage{adjustbox}
\usepackage{epstopdf}
\usepackage[acronym,toc,shortcuts]{glossaries}
\usepackage{cite}
\usepackage{rotating}
\usepackage{eurosym}
\usepackage{amsfonts}
\usepackage{mdwmath}
\usepackage{mdwtab}
\usepackage{eqparbox}
\usepackage{float}
\usepackage[utf8]{inputenc}
\usepackage{tabularx}
\usepackage{txfonts}

\input{acronyms.tex}
\makeglossaries

\begin{document}

\title{Enhancing logical deduction with math: the rationale behind Gardner and Carroll}

\author{\IEEEauthorblockN{Tiziana Castellano\IEEEauthorrefmark{1}, Pietro Boccadoro\IEEEauthorrefmark{3}\IEEEauthorrefmark{2}}
\IEEEauthorblockA{\IEEEauthorrefmark{3}Dep. of Electrical and Information Engineering (DEI), Politecnico di Bari, Bari, Italy}
Email:name.surname@poliba.it\\
\IEEEauthorrefmark{2}CNIT, Consorzio Nazionale Interunivesitario per le Telecomunicazioni, Politecnico di Bari, Italy\\
\IEEEauthorrefmark{1}Email:surname.name@gmail.com}

\maketitle

\begin{abstract}
Math is widely considered as a powerful tool and its strong appeal depends on the high level of abstraction it allows in modelling a huge number of heterogeneous phenomena and problems, spanning from the static of buildings to the flight of swarms. As a further proof, Gardner’s and Carroll’s problems have been intensively employed in a number of selection methods and job interviews. Despite the mathematical background, these problems are based on, several solutions and explanations are given in a trivial way. This work proposes a thorough investigation of this framework, as a whole. The results of such study are three mathematical formulations that express the understood mathematical relationship in these well-known riddles. The proposed formulas are of help in the formalization of the solutions, which have been proven to be less time-taking when compared to the well-known classic ones, that look more heuristic than rigorous.
\end{abstract}


\IEEEpeerreviewmaketitle

\section{Introduction}\label{sec:intro}
Mathematics is widely considered as a powerful tool. Its appeal strongly depends on its ability to clearly describe both abstract and high-level models and real-world phenomena. For instance, it is used to efficiently describe heterogeneous problems, spanning from the static of buildings to the flight of birds.
This flexible subject has so many application that also include puzzles, riddles and a number of different problems.
Among the many applications math is employed into there are also selection tests, often based on puzzles and problems that can be solved relying on algebra, trigonometry and probability theory.
This kind of problems is so widespread that has become part of calls for selection and job interviews, specifically aimed at verifying problem solving attitude rather than the ability to lateral thinking.
These selection criteria are mainly based on what has been studied and proposed by Lewis Carroll and Martin Gardner \cite{carroll1958mathematical,gardner1985misteri,gardner1986entertaining}. The books were written to propose some mind blowing puzzles that could only be solved relying on a mathematical background. 

In this work, a detailed study of some of those problems is proposed. The "classic" solutions for each are herein compared to the one that can be obtained using the three proposed equations, i.e., Titti's formulas, which are specifically conceived to address the solution of three kind of problems, for instance, (i) \textit{three-based problem}, (ii) \textit{heavy-weight problem}, and (iii) \textit{socks problem}.
For each one, the "classic" solution, i.e., the one that can be found in books used to study for tests and competitions \cite{artquiz18, artquiz182,unidtest1,unidtest2,editest18,hoepli16}, is discussed.
The proposed work demonstrates how the simple application of the formulas can ease the resolution of the same riddles, thus making it faster while increasing the formal mathematical formulation. 

The remainder of this work is as follows: Section \ref{sec:background} presents the background and related works. Section \ref{sec:problems} reports several problems and puzzles, divided into three categories, together with their solution, explained in details. In Section \ref{sec:proposal}, instead, the same problems are solved using the proposed formulas. The mathematical formulations presented in this contribution are applied to the same set of problems to demonstrate their effectiveness. Finally, Section \ref{sec:conclusions} concludes the work and draws future possibilities.

\section{Background and Motivation}\label{sec:background}
Math has always been a tricky but fascinating topic to study and puzzles can be considered as a clear example of this matter of fact.
Among the many that experimented such fascination, Charles L. Dodgson (who is better known by the pseudonym of Lewis Carroll) is probably one of the brightest example. It has been a while since its book "Pillow-Problems" has been published \cite{carroll1958mathematical}.
This publication contains more than seventy mathematical puzzles, pretty original at the time.
While Carroll had worked all of them out in his head and encouraged the reader to do the same, many of the puzzles are not easy at all.
The problems span from those that can be solved with simple algebra or plane geometry, to those that require more advanced math concepts, i.e., trigonometry and probability.
Among those related to probability, one of the most popular involves some containers and some objects with different colors. In this kind of problem, some of the objects with a certain color are moved from one container to another. The problem asks to calculate the probability of drawing one of the objects moved to the new container.
The general answer to this type of problem is $2n / (n + d)$, where n in the number of objects with the same color within one of the two containers, whereas d is the total amount of colored objects in the other one.
This formula can be considered as an equivalent of the Bayes' Theorem.


In \cite{gardner1985misteri}, Martin Gardner explains mathematical puzzles that have been published on several journals, such as Scientific American and Science World Journal.
The author discusses the magic tricks from a mathematical point of view. A whole chapter is dedicated to the geometric disappearances treated by a very logical point of view.

In \cite{dudeney1951amus}, Dudeney proposes a number of puzzles that the author claim could be solved using logical reasoning and maths.
The reference puzzles do not respect formalism. Still they are of reference as the simplest ones are tricky whereas the more complex could be discussed in mathematical terms. Most of the puzzles are discussed in terms of arithmetic and geometric.

Mathematical puzzles is a set of both old and new puzzles, grouped into sections, covering a variety of mathematical topics: plane and solid geometry, probability, topology \cite{gardner1986entertaining}.
For instance, the probability section emphasizes the role and importance of probability laws. Geometric puzzles, instead, are proposed as to test the ability to think in more than one dimension. As for the topology, among the "youngest and most disordered branches of modern geometry", offers a look at a strange dimension in which the properties remain unchanged, regardless of how a figure is twisted, stretched or compressed.
The book suggests that almost everything that is needed to solve a riddle is the ability to think logically and clearly.
Among the most popular, there is one about distance and time, where there is a subject who arrives at a given time at the station, and another person who goes to pick him up by car, to go together to a further destination. One day, the first subject arrives as first at the station and sets out alone. The two arrive at their destination in advance. The request is related to the time in which the subject walks alone before meeting the second person. According to Gardner, the answer to this question is $X-(Y/2)$, where X indicates how much time before the man has arrived at the station before the first subject while Y indicates how far ahead they reach their destination.

Both \cite{gardner1992fractal} and \cite{gardner2015giochi} are collections of extracts from Gardner's' "Scientific American" column. Each article deals with new ideas and new solutions regarding its puzzles. The highlights include two new chapters, on poetry and minimal sculpture.

\cite{seneta1993lewis} is an article, written by Seneta, which only deals with one of the probability themes studied by Carroll in "Pillow Problems", or on binomial probability distributions to determine, for example, the bearing of extracting a colored object from a container, previously moved in it from another container referred to in \cite{carroll1958mathematical}.
On the same topic, Vakhania wrote \cite{vakhania2009probability}, a publication on similar themes in which he deduced a general mathematical formulation.

\cite{geronimi2017alice} is a book dedicated to Gardner. Nando Geronimi restate the importance of Gardner's "discovery" for his education as a teacher and a great lover of mathematical games. Some of the games conceived or promoted by Gardner have become true classics: Maurizio Paolini comments on the "tumbling rings", Alessandro Musesti instead speaks of "Life" (born from an idea of the mathematician John Conway).

\section{The problems}\label{sec:problems}
This section identifies three different kind of puzzles and describes the way the solutions are generally proposed.
\subsection{3-based puzzle}
In the 3-based puzzles, the involved entities (i.e., the variables) are in a number of three and the key to solve the puzzle lies beyond the relationship between them.
For example: 
\begin{itemize}
	\item \textit{Problem 1}:
	If six (6) cats eat six (6) mice in six (6) minutes, how many cats eat one hundred (100) mice in fifty (50) minutes?
		
	\item \textit{Problem 2}:
	If one hundred (100) cats eat one hundred and fifty (150) mice in one (1) hour, home many cats eat sixty (60) mice in thirty (30) minutes?
	
	\item \textit{Problem 3}:
	If three (3) bakers bring out forty (40) sandwiches in one hundred and twenty (120) minutes, how many bakers will bring out one hundred (100) sandwiches in thirty (30) minutes?
\end{itemize}

The solutions to the proposed puzzles are as follows:
\begin{itemize}
	\item \textit{Solution 1}: Given that 6 cats eat 1 mouse every minute, it means 2 mice every 2 minutes, 50 mice in 50 minutes and so on. Therefore, it takes 12 cats to eat 100 mice in 50 minute.
	
	\item \textit{Solution 2}: If 100 cats eat 150 mice in 1 hour, it means that 100 cats will eat 75 mice in half an hour. In case mice were 60, it clearly results that: 
	$ 100 / 75 = x / 60 $ which implies $ x = (6000)/75 = 80 $.
	Therefore, 80 cats eat 60 mice in 30 minutes.
	
	\item \textit{Solution 3}: A baker is able to take out from the oven 40/3 sandwiches in 120 minutes\footnote{Reference work quantity.}. Therefore, the baker is able to take out  $ 40/3 * 1/120 = 1/9 $ sandwiches every minute\footnote{Reference work quantity in a minute}. Therefore, if a baker takes out 1/9 sandwiches every minute, in 30 minutes he will take out $ 1/9 * 30= 30/9 = 10/3 $ sandwiches, which leads to $ 10/3 : 1 = 100 : x $. Therefore, the solution for the problem is:
	$ x= 1*100/(10/3) = 100*3/10=30 $.	
\end{itemize}

\subsection{Heavy weight puzzle}
In this kind of problems, it is generally given a group of objects that are all equal one with the other except for one that is heavier than the others.
The name of this puzzle clearly expresses the understood question: \textit{how could you identify the odd-one-out (i.e., the heavier)?} This question is often formulated as: \textit{how many times do you need to weight those objects before you identify the heavier?}
Nevertheless, it is not always that simple, since sometimes the question changes to \textit{which is the lowest number of weights you need to identify the heavier?}
Here, three problems of this kind are proposed, first, and solved, afterward:
\begin{itemize}
	\item \textit{Problem 4}:
	Given thirteen (13) coins, twelve (12) coins are equal and only one (1) of them is different in weight (i.e., heavier). How many weights it take to find the odd-one-out?
	
	\item \textit{Problem 5}:
	Lorna has five (5) necklaces. One out of the five is heavier than the others. Having a two-plate scale, how many weights it takes to identify the odd-one-out?
	
	\item \textit{Problem 6}:
	Nelly has nine (9) stamps, apparently equal one with the other. Nevertheless, one of them is heavier than the others. With a two-plate scale, how many weighs will it take to find out the heavier one?
\end{itemize}

The solutions for the three are as follows:
\begin{itemize}
	\item \textit{Solution 4}:
	The first step toward the solution consists of an equal subdivision of the weights into groups. In case the objects are odd, there will be only one group with an extra object. If they are even, there will be a group with one object less. In this problem, the objects are thirteen, therefore there will be three different groups, two composed by 4 elements and one counting 5 elements\footnote{If the objects were fourteen, the three groups would have been composed by five, five and four elements, respectively.}.
	The first weighing between two out of three groups is done, while another one is set aside. In general, the group taken aside is the one with the extra token.
	In this phase, two possible outcomes may occur: in the former, the weighing reveals that the groups are different, therefore the odd-one-out is within the heavier group. In the latter, instead, the two are equal.
	In any case, the outcome will spot out which group contain the heavier coin and, from that point on, the second subdivision will be carried out on that group.
	In the former case, the subsequent subdivisions will split the group into 2 groups of 2 coins and 2 groups with only 1 coin each. This weighing will be the last one and will surely spot out the odd-one-out, thus leading to the correct solution in a total of 3 weights.
	In the latter case, instead, the group to be split will be the one counting 5 elements. In this case, the subdivision will create 2 groups made by 2 coins and one group counting only 1 element. The aforedescribed procedure is repeated once again for the groups of 2.
	It clearly comes out that no matter the group the heavier coin is within, the maximum number of weighing\footnote{In case they were fourteen, or even, the process is repeated in the very same way.} will be 3.

	\item \textit{Solution 5}:
	The procedure that leads to the solution of this problem is pretty much the same as the previous one, starting from the second subdivision. It clearly results that the needed weighings are 2.
	
	\item \textit{Solution 6}:
	To solve the problem, the nine stamps are equally divided into three groups of three elements. In this case, by picking any two out of the three groups, the weighing will clearly spot out the one that contains the heavier\footnote{It i worth noting that the very same thing happens for any group with a numerosity that results to be a power of 3.}. Therefore, the method will continue subdividing the identified group of interest into three groups of 1 element. Whichever the outcome, the heavier stamp will be spotted out.
\end{itemize}

\begin{figure*}[htbp]
	\centering
	\includegraphics[width=0.6\textwidth]{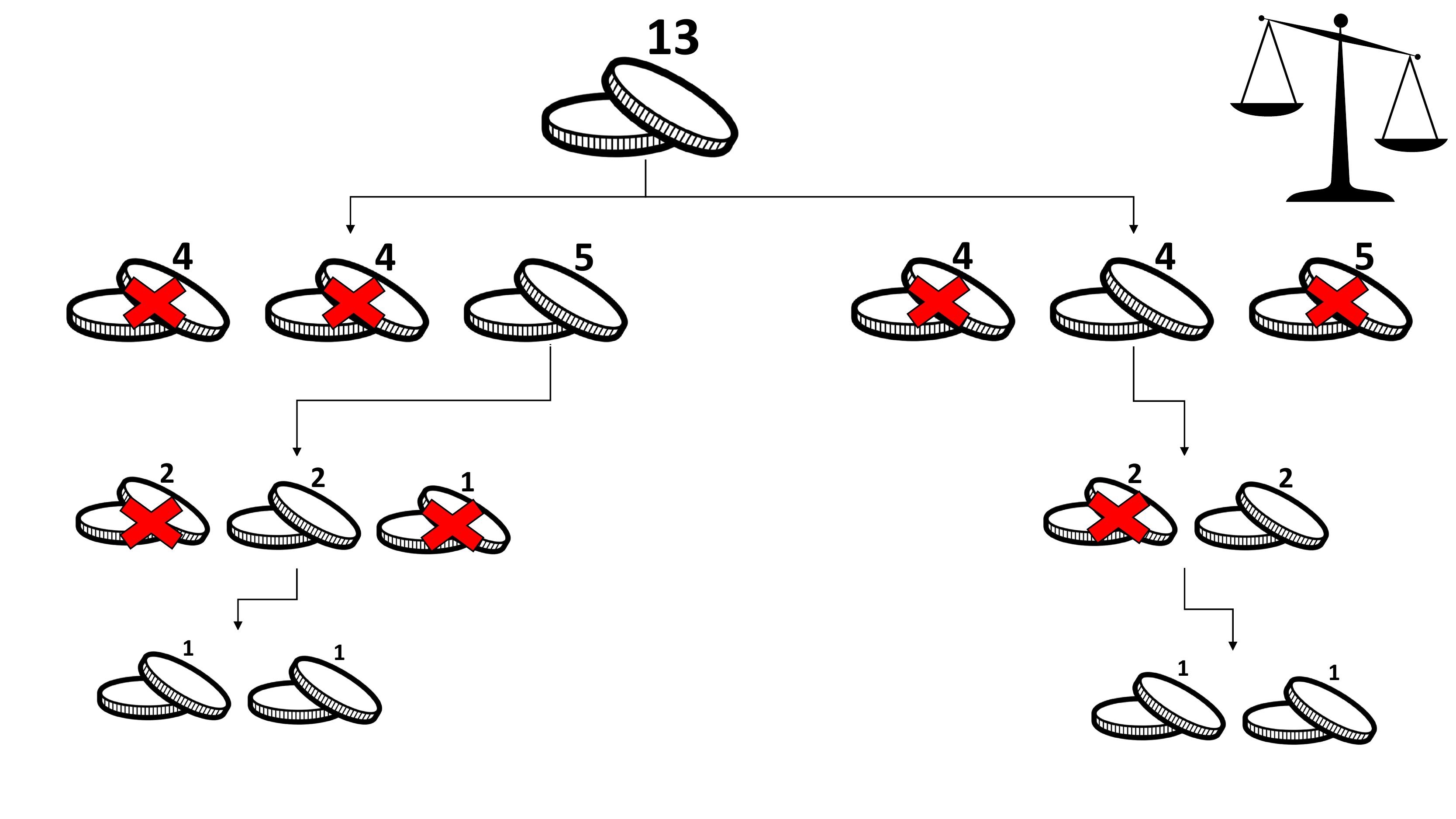}
	\caption{Graphical representation of the solution for the heavy-weight problem.}
	\label{fig:scale}
\end{figure*}

In Figure \ref{fig:scale}, an abstract representation of the way Heavy weight puzzles can be solved is given.

\subsection{Socks puzzle}
The third category of puzzles of interest in this work is the so-called \textit{socks puzzles}, that is generally proposed as follows:
\begin{itemize}
	\item \textit{Problem 7}:
	In a drawer, there are five (5) pairs of blue socks, four (4) pairs of red socks and six (6) pairs of black socks. Randomly pulling out socks, how many catches will it take to have a pair of the same color?
	
	\item \textit{Problem 8}:
	A tailor, late in delivering a fancy jacket for a customer, needs a set of four (4) identical buttons. He instructs his son to bring him four (4) loaves of the same color chosen in a drawer containing eighty-four (84) blue buttons, thirty-two (32) turquoise, twenty-eight (28) red and four (4) green, all of the same shape and size. Since tailor's son cannot distinguish the colors in the dark, which is the minimum number of buttons the tailor will he have to take to have four (4) buttons of the same color?
\end{itemize}

The general solutions to those problems are as follows:
\begin{itemize}
	\item \textit{Solution 7}:
	The first sock extracted will be of a certain color. Once the second is extracted, it might be of the same color or not. The former case implies a trivial solution. In the latter, instead, the sock will be of a different color. When it comes to the third extraction the non-trivial possibility implies that the this third sock is of a third color too. In this case, no matter which color is extracted, it will certainly be equal to one of the previously taken out.
	As an example, a possible solution to the problem could be: blue - red - black - blue.
	
	\item \textit{Solution 8}:
	The solution to this problem can be described as previously done. Given a sequence of four extraction of four buttons of four different colors, the fifth extracted button will certainly be of a previously extracted one. In case the sequence repeats itself, the correct number is 13.
	As an example, a possible solution is: blue - turquoise - red - green - blue - turquoise - red - green - blue - turquoise - red - green - blue.
\end{itemize}

\section{The proposed approaches}\label{sec:proposal}
The solutions seen so far are surely leading to answering the proposed questions. Nevertheless, they do not look formal as math should be. Moreover, despite the correctness of the answer, the way the solutions are given could be missleading and/or misunderstood. Furthermore, they do not seem to explain the mathematical relationships that the puzzles are based onto.
The present contribution proposes a mathematical formulation for the three set of puzzles described in Section \ref{sec:problems}.
With Table \ref{tab:symbols}, the main symbols used herein are summarized.
\begin{table}[htbp]
	\centering
	\begin{tabular}{|c|l|}
		\hline
		\textbf{Symbol} & \multicolumn{1}{c|}{\textbf{Explanation}} \\ \hline
		W               & Work                                      \\ \hline
		t               & Time                                      \\ \hline
		S               & Subjects                                  \\ \hline
		N               & Number of objects                         \\ \hline
		i               & Exponent                                  \\ \hline
		P               & Number of weight counts                   \\ \hline
		$n_C$            & Number of available colors                           \\ \hline
		$n_R$            & Number of requested socks                        \\ \hline
	\end{tabular}
	\caption{Summary of symbols and notation.}
	\label{tab:symbols}
\end{table}

\subsection{First Titti's Formula}
The first formula solves the \textit{3-based problems}:
\begin{equation} \frac{W}{S t} = k \label{formula:1}\end{equation}
being k a constant value.
The three reference quantities are independent and proportional one to the other.
The approach to the solution can be similar to what happens with physics problems, that can be solved varying one out of the three while calculating the resulting values for the others.

For example, given the work (W), both the subjects (S) and the time (t) will vary. The relationship between the latter two is an inverse proportionality, which means that as long and the time rises, a lower number of subjects is needed for the work to be completed. Counterwise, the more subjects are involved, the less time it will be necessary for the work to be completed. Such relationship lead to:

\begin{equation} S t = constant \label{eq:2} \end{equation}

In a similar fashion, given the time, it is possible to evaluate the relationship between the work and the involved subjects. It can be easily understood that the higher the number of subjects involved, the more the produced work. Therefore, the number of subjects and the time taken are directly proportional. 

\begin{equation} W/S = constant\label{eq:3} \end{equation} 

Lastly, if the number of involved subjects is given, both work and time can be estimated. The way the problem is solved is pretty similar to the previous cases. For instance, if time rises, so will the work, and viceversa:

\begin{equation} W/t = constant \end{equation} 

By jointly considering (\ref{formula:1}), (\ref{eq:2}) and (\ref{eq:3}), it results that:

\begin{equation} W/(S t) = constant \end{equation}

Therefore, given these three values, it is quicker to solve the problems.
Turning from theory to practice, all the problems proposed and solved in Section \ref{sec:problems} are herein solved once again here, using the proposed formulation:

\begin{itemize}
	\item \textit{Solution 1 with the formula}:
	Given the number of cats (6), minutes (6) and mice (6), the variable is the number of cats that are able to eat 100 mice in 50 minutes.
	With Formula (\ref{formula:1}), it can be assumed that work (W) = mice, subjects (S) = cats. Therefore:
	$ 6/(6*6) = 100/(50*x) $ which results in $ x = 12 $.
	
	\item \textit{Solution 2 with the formula}:
	Given the number of cats(100), minutes (60), and mice (150), the variable is the number of cats that eat 60 mice in 30 minutes.
	Using Formula (\ref{formula:1}), it results that:
	$ 150/(100*60) = 60/(30*x)$, therefore $ x = 80 $.
	
	\item \textit{Solution 3 with the formula}:
	In this third case, it is possible to assume that:
	$W = sandwiches = 40$, $S = bakers = 3$, and $t = 120 min$.
	Being the number of bakers the variable of interest, thanks to Formula (\ref{formula:1}), it can be obtained $ 40/(3*120) = 100/(30*x) $ from which it is obtained $ x = 30 $.
\end{itemize}

\subsection{Second Titti's Formula}
The second formula solves the \textit{heavy weight puzzle} and is as follows:
\begin{equation} 3^i < N \leq 3^{i+1} \label{formula:2} \end{equation}
with $i \in \mathbb{N}$.
The correct answer to the question is:
$ P = i+1 $

Also in this case, all the problems proposed and solved in Section \ref{sec:problems} are solved again here, using the proposed formulation:
\begin{itemize}
	\item \textit{Solution 4 with formula}:
	The number of objects N is equal to 13. Using Formula (\ref{formula:2}), it is possible to write $ 3^2<13\leq 3^3$. Therefore, the solution is $P=3$.
	
	\item \textit{Solution 5 with formula}:
	In this case, the number of objects is equal to 9, which can also be expressed as $3^2$. Using Formula (\ref{formula:2}), it is possible to write $ 3^1<9\leq 3^2$ which leads to the solution $P=2$.
	
	\item \textit{Solution 6 with formula}:
	In this problem, N is equal to 4. Using Formula (\ref{formula:2}), it is possible to write $ 3^1<4\leq 3^2$, from which it results that $P=2$.
\end{itemize}
It is worth noting that the application of this formula can be considered as a noticeable aid in shortening the time taken for solving the problems. This is even more evident when the objects are more than 9.

\subsection{Third Titti's Formula}
The third formula proposed in this contribution is:
\begin{equation} [n_C (n_R-1)]+1 \label{formula:3}\end{equation}
and solves the \textit{socks puzzle}.

Again, all the problems proposed and solved in Section \ref{sec:problems} are solved here using the proposed formulation:
\begin{itemize}
	\item \textit{Solution 7 with formula}:
	The number of colors is $ n_C = 3$, the number of required socks of the same color is equal to 2. With Formula (\ref{formula:3}) it results: $[3(2-1)]+1 = 4$.
	
	\item \textit{Solution 8 with formula}:
	The number of different colors is equal to 4, whereas the number of required buttons of the same color is 4. When applying Formula (\ref{formula:3}), it results that: $[4(4-1)]+1 = 13$.
\end{itemize}

It is worth specifying that the proposed formula is of great relevance when the number of indicated objects of the same color is greater than three. In this case, in fact, the aforedescribed sequences are not needed, thus resulting in a quicker solution for the problem.

\section{Conclusions}\label{sec:conclusions}
This work presented three mathematical formulations to express the understood mathematical relationship between the entities involved in some of the most well-known Gardner's and Carroll's problems. The formulas can be considered of help in the formalization of the solutions for a number of reasons. First of all, the solution do no longer require dedicated graphs and/or a written or understood reasoning, as it can be solved by simply applying the proper formula. Moreover, these formulations have been proven to be less time-taking when compared to the well-known classic ones, that look more heuristic than rigorous. 
\bibliography{Bibliography}
\bibliographystyle{IEEEtran}

\end{document}

%% file: acronyms.tex
\newacronym{6lo}{6lo}{\ac{IPv6} over networks of resource-constrained nodes}
\newacronym{6LoWPAN}{6LoWPAN}{\ac{IPv6} over Low power Wireless Personal Area Networks}
\newacronym{6top}{6top}{6tisch Operation Sublayer}
\newacronym{6TiSCH}{6TiSCH}{IPv6 over the TSCH mode of IEEE 802.15.4}
\newacronym{AGV}{AGV}{Automated Guided Vehicle}
\newacronym{AP}{AP}{Allocation Policy} 
\newacronym{AQ}{AQ}{Autonomous Quadcopter}
\newacronym{ASAP}{ASAP}{As-Soon-As-Possible Schedule Another Parent}
\newacronym{ASN}{ASN}{Absolute Slot Number}
\newacronym{ASV}{ASV}{Autonomous Surface Vehicle}
\newacronym{AUV}{AUV}{Autonomous Underwater Vehicle}
\newacronym{CCTV}{CCTV}{Closed Circuit TeleVision}
\newacronym{CEA}{CEA}{Cell Estimation Algorithm}
\newacronym{CIR}{CIR}{Channel Impulse Response}
\newacronym{CoAP}{CoAP}{Constrained Application Protocol}
\newacronym{CoRE}{CoRE}{Constrained RESTful Environments}
\newacronym{CR}{CR}{Cognitive Radio}
\newacronym{CSMA/CA}{CSMA/CA}{Carrier Sense Multiple Access with Collision Avoidance}
\newacronym{CTD}{CTD}{Conductivity, Temperature, Depth}
\newacronym{D2D}{D2D}{Device-to-Device}
\newacronym{DAG}{DAG}{Directed Acyclic Graph}
\newacronym{DAO}{DAO}{Destination Advertisement Object}
\newacronym{DETB}{DETB}{Dynamic Energy \& Traffic Balance}
\newacronym{DGAC}{DGAC}{Direction Générale de l'Aviation Civile}
\newacronym{DIO}{DIO}{Destination Information Object}
\newacronym{DNS}{DNS}{Domain Name System}
\newacronym{DODAG}{DODAG}{Destination Oriented DAG}
\newacronym{DSSS}{DSSS}{Direct Sequence Spread Spectrum}
\newacronym{E2E}{E2E}{End-to-End}
\newacronym{EASA}{EASA}{European Aviation Safety Agency}
\newacronym{ENAC}{ENAC}{Ente Nazionale per l'Aviazione Civile}
\newacronym{EOF}{EOF}{End Of Frame}
\newacronym{ePFC}{ePFC}{extended Potential Field Controller}
\newacronym{FAA}{FAA}{Federal Aviation Administration}
\newacronym{FANET}{FANET}{Flying Ad-hoc NETwork}
\newacronym{FFD}{FFD}{Fully Functioned Device}
\newacronym{FSO}{FSO}{Free Space Optical}
\newacronym{GIS}{GIS}{Geographic Information System}
\newacronym{GPS}{GPS}{Global Position System}
\newacronym{HTTP}{HTTP}{HyperText Transfer Protocol}
\newacronym{ICN}{ICN}{Information-Centric Networking}
\newacronym{ICT}{ICT}{Information and Communication Technologies}
\newacronym{IETF}{IETF}{Internet Engineering Task Force}
\newacronym{IIoT}{IIoT}{Industrial Internet of Things}
\newacronym{IMC}{IMC}{Internal Model Control}
\newacronym{IMU}{IMU}{Inertial Measurement Unit}
\newacronym{IoD}{IoD}{Internet of Drones}
\newacronym{IoT}{IoT}{Internet of Things}
\newacronym{IP}{IP}{Internet Protocol}
\newacronym{IPv6}{IPv6}{Internet Protocol version 6}
\newacronym{ISM}{ISM}{Industrial, Scientific and Medical}
\newacronym{ITS}{ITS}{Intelligent Transportation System}
\newacronym{JSON}{JSON}{JavaScript Object Notation}
\newacronym{LLN}{LLN}{Low-power Lossy Network}
\newacronym{LLSF}{LLSF}{Low-power Lossy Network}
\newacronym{LoRaWAN}{LoRaWAN}{Long Range Wide Area Network}
\newacronym{LoS}{LoS}{Line of Sight}
\newacronym{LPWAN}{LPWAN}{Low Power Wide Area Network}
\newacronym{LTE}{LTE}{Long-Term Evolution}
\newacronym{M2M}{M2M}{Machine-to-Machine}
\newacronym{MAC}{MAC}{Medium Access Control}
\newacronym{MANET}{MANET}{Mobile Ad-hoc NETwork}
\newacronym{MCU}{MCU}{Micro-Controller Unit}
\newacronym{MSPRT}{MSPRT}{Multi Sequential Probability Ratio Test}
\newacronym{MTC}{MTC}{Machine-Type Communication}
\newacronym{mcMTC}{MTC}{mission critical Machine-Type Communication}
\newacronym{MTU}{MTU}{Maximum Transmission Unit}
\newacronym{NLoS}{NLoS}{Non-Line of Sight}
\newacronym{O-QPSK}{O-QPSK}{Offset-Quadrature Phase-Shift Keying}
\newacronym{OS}{OS}{Operating System}
\newacronym{OTF}{OTF}{On The Fly}
\newacronym{P2P}{P2P}{Peer-to-Peer}
\newacronym{PCC}{PCC}{Path Computation Client}
\newacronym{PCE}{PCE}{Path Computation Element}
\newacronym{PCEP}{PCEP}{Path Computation Element Protocol}
\newacronym{PDR}{PDR}{Packet Delivery Ratio}
\newacronym{PHY}{PHY}{Physical Layer}
\newacronym{PID}{PID}{Proportional-Integral-Derivative}
\newacronym{PLA}{PLA}{PolyLactic Acid}
\newacronym{PLR}{PLR}{Packet Loss Ratio}
\newacronym{PxIMU}{PxIMU}{Pixhawk Inertial Measurement Unit}
\newacronym{RFD}{RFD}{Reduced Functioned Device}
\newacronym{RMM}{RMM}{Random Mobility Model}
\newacronym{ROLL}{ROLL}{Routing Over Low power and Lossy networks}
\newacronym{ROS}{ROS}{Robot Operating System}
\newacronym{ROV}{ROV}{Remotely Operated underwater Vehicle}
\newacronym{RPL}{RPL}{Routing Protocol for Low-power and Lossy networks}
\newacronym{RTT}{RTT}{Round Trip Time}
\newacronym{RSSI}{RSSI}{Received Signal Strength Indication}
\newacronym{QoE}{QoE}{Quality of Experience}
\newacronym{QoL}{QoL}{Quality of Link}
\newacronym{QoS}{QoS}{Quality of Service}
\newacronym{SAIN}{SAIN}{Space-Air Integrated Network}
\newacronym{SDN}{SDN}{Software Defined Network}
\newacronym{SFD}{SFD}{Start of Frame Delimiter}
\newacronym{SF0}{SF0}{Scheduling Function Zero}
\newacronym{SINR}{SINR}{Signal to Interference-plus-Noise Ratio}
\newacronym{SLAM}{SLAM}{Simultaneous Localization And Mapping}
\newacronym{SNR}{SNR}{Signal-to-Noise Ratio}
\newacronym{SPI}{SPI}{Serial Peripheral Interface}
\newacronym{TCP}{TCP}{Transmission Control Protocol}
\newacronym{TDMA}{TDMA}{Time Division Multiple Access}
\newacronym{TI}{TI}{Texas Instruments}
\newacronym{TSCH}{TSCH}{Time Slotted Channel Hopping}
\newacronym{UDP}{UDP}{User Datagram Protocol}
\newacronym{UASN}{UASN}{Underwater Acoustic Sensor Network}
\newacronym{UAV}{UAV}{Unmanned Aerial Vehicle}
\newacronym{UGV}{UGV}{Unmanned Ground Vehicle}
\newacronym{UWSN}{UWSN}{Underwater Wireless Sensor Network}
\newacronym{USV}{USV}{Unmanned Surface Vehicle}
\newacronym{UUV}{UUV}{Unmanned Underwater Vehicle}
\newacronym{UTM}{UTM}{Unmanned Aerial System Traffic Management}
\newacronym{V2I}{V2I}{Vehicle-to-Infrastructure}
\newacronym{V2V}{V2V}{Vehicle-to-Vehicle}
\newacronym{VANET}{VANET}{Vehicular Ad-hoc NETwork}
\newacronym{VLC}{VLC}{Visible Light Communication}
\newacronym{WG}{WG}{Working Group}
\newacronym{WSN}{WSN}{Wireless Sensor Network}
\newacronym{CH}{CH}{Cluster Head}